\documentclass{amsproc}

\usepackage{amssymb,amsmath,amsfonts,amsthm}
\newtheorem*{thm}{Theorem}
\newtheorem*{lem}{Lemma}

\newtheorem*{ac}{Acknowledgements}
\theoremstyle{remark}
\newtheorem*{remark}{Remark}
\theoremstyle{remark}

\newcommand{\ext}{\mathop{\mathrm{Ext}}\nolimits}

\newcommand{\onto}{\twoheadrightarrow}

\newcommand{\iso}{\simeq}

\newcommand{\projdim}{\mathop{\mathrm{proj\thinspace dim}}\nolimits}
\newcommand{\injdim}{\mathop{\mathrm{inj\thinspace dim}}\nolimits}
\newcommand{\tor}{\mathop{\mathrm{Tor}}\nolimits}

\newcommand{\x}{\mathbf x}
\newcommand{\longto}{\longrightarrow}

\renewcommand{\hom}{\mathop{\mathrm{Hom}}\nolimits}

\renewcommand{\tilde}{\protect\widetilde}

\renewcommand{\hat}{\protect\widehat}

\begin{document}
\title[CHARACTERIZATION OF MODULES OF FINITE PROJECTIVE DIMENSION]
{CHARACTERIZATIONS OF MODULES OF FINITE PROJECTIVE DIMENSION OVER
COMPLETE INTERSECTIONS}
\author{Jinjia Li}
\address{Department of Mathematics, University of Illinois at Urbana-Champaign, Urbana, Illinois 61801}
\curraddr{Department of Mathematics, University of Illinois at
Urbana-Champaign, 1409 W. Green Street, Urbana, Illinois 61801}
\email{jinjiali@math.uiuc.edu}
\thanks{This research was carried out while the author was supported by a research grant from the
UIUC Campus Research Board of the University of Illinois under the supervision of S. Dutta.}
\subjclass{Primary 13C14, 13C40, 13D05, 13D40, 13H10.}

\date{July 1, 2004}

\keywords{complete intersection, finite projective dimension,
finite injective dimension, flatness, Frobenius, Ext, Tor.}

\begin{abstract}
Let $M$ be a finitely generated module over a local complete
intersection $R$ of characteristic $p>0$. The property that $M$
has finite projective dimension can be characterized by the
vanishing of $\ext_R^i({}^{f^n}\!\! R,M)$ for some $i>0$ and for
some $n>0$.
\end{abstract}

\maketitle Let $(R,m,k)$ be a local ring of characteristic $p>0$.
The Frobenius endomorphism $f_R:R \to R$ is defined by
$f_R(r)=r^p$ for $r \in R$. Each iteration $f_R^n$ defines a new
$R$-module structure on $R$, denoted by ${}^{f^n}\!\! R$ for which
$a\cdot b=a^{p^n}b$. For any $R$-module $M$, $F_R^n(M)$ will stand
for $M\otimes_R {}^{f^n}\!\! R$ and $\tilde F_R^n(M)$ will stand
for $\hom_R({}^{f^n}\!\! R,M)$. Avramov and Miller [A-M] proved
that over a local complete intersection ring, if a finitely
generated $R$-module $M$ satisfies $\tor_i^R(M,{}^{f^n}\!\! R)=0$
for some fixed $i,n>0$, then it is of finite projective dimension.
Later, Dutta [D] provided a simple proof of this result without
using the notion of ``complexity''. Using a similar method as in
[D], we obtain another characterization of finitely generated
modules of finite projective dimension over complete intersection
rings.

\medskip

\begin{thm}Let $M$ be a finitely generated module over a local complete
intersection ring $R$. If for some $i,n>0$, $\ext_R^i({}^{f^n}\!\!
R,M)=0$, then $M$ has finite injective dimension over $R$.
\end{thm}

\begin{proof}Without loss of generality we can assume that
$R$ is complete and the residue field $k$ is perfect. Let $R=S/\x
$ where $S$ is a complete regular local ring of characteristic
$p>0$ and $\x=(x_1,\dots,x_r)$ is an ideal generated by an
$S$-sequence $x_1,\dots,x_r$. Write $R_n=S/\x^{p^n}$. Notice that
with this notation, every $R_n$-module is also an $R_l$-module via
the natural surjection $R_l \onto R_n$ for all $l>n$.

Let $\tilde{f}^n=f^n_S \otimes_S R$ (base change of $f^n_S: S \to
S$ along the $S$-algebra $R$). We know by Kunz's Theorem [K,
Theorem 3.3] that $f_S^n$ is flat. Since $k$ is perfect, $f_S^n$
is module finite. It follows that $\tilde{f}^n: S/\x \to
S/\x^{p^n}$ is also flat and module finite. Observe that the map $
R \overset{f_R^n}{\longto} R$ can be factored as
\[
S/\x \overset {\tilde{f}^n}{\longto} S/\x^{p^n} \overset
{\eta_n}{\longto} S/\x, \tag{1}
\]
or
\[R \overset {\tilde{f}^n}{\longto} R_n \overset
{\eta_n}{\longto} R \] where $\eta_n: S/\x^{p^n} \longto S/\x$ is
the natural surjection.

Thus $f^n_R=\eta_n \cdot \tilde {f}^n$.

It follows from the adjointness of $\hom$ and Tensor that for any
$R$-module $T$,
\[\hom_{R_n}(R,\hom_R({}^{\tilde{f}^n}\!\!
R_n,T)) \iso \hom_R(R \otimes_{R_n}{}^{\tilde{f}^n}\!\! R_n,T)
\iso \hom_R({}^{f^n}\!\! R,T). \tag{2}\] Let $I^\bullet$ be an
injective resolution of $M$ over $R$. Since $\tilde {f}^n$ is
flat, $\hom_R({}^{\tilde{f}^n}\!\! R_n,I^\bullet)$ is an injective
resolution of $\hom_R({}^{\tilde{f}^n}\!\! R_n,M)$ over $R_n$.
Therefore, by (2)
\begin{align*}
\ext_R^j({}^{f^n}\!\! R, M)
&=H^j(\hom_R({}^{f^n}\!\! R,I^\bullet))\\
&\iso H^j(\hom_{R_n}(R,\hom_R({}^{\tilde{f}^n}\!\! R_n,I^\bullet)))\\
&=\ext_{R_n}^j(R,\hom_R({}^{\tilde{f}^n}\!\! R_n,M)), \forall j
\geq 0.
\end{align*}
But
\begin{align*}
\hom_R({}^{\tilde{f}^n}\!\! R_n,M)&=\hom_R({}^{f^n}\!\!
S\otimes_s{R},M)\\
                                  &\iso \hom_S({}^{f^n}\!\! S,M)\\
                                  &= \tilde F_S^n(M).
\end{align*}
Hence, $\ext_R^j({}^{f^n}\!\! R,M)\iso \ext _{R_n}^j(S/\x, \tilde
F_S^n(M)), \forall j \geq 0$.

\medskip

Next, we show $\ext_R^{i}({}^{f^n}\!\! R,M)=0$ for some $i>0$
implies $\ext_R^{i}({}^{f^n}\!\! R,M)=0$ for all $i>0$. Consider a
filtration of $S/\x^{p^n}$ of the form:
\[
0\to K_1\to S/\x^{p^n}\to S/\x\to 0,
\]
\[
0\to K_2\to K_1\to S/\x\to 0,
\]
\[
\vdots
\]
\[
0\to K_{t_n}\to K_{t_{n-1}}\to S/\x\to 0,
\]
where $K_{t_n}=S/\x$. By applying $\hom_{R_n}(-,\tilde F_S^n(M))$
to the above short exact sequences, we get the long exact
sequences of $\ext$'s. Since $\ext_{R_n}^i(S/\x,\tilde
F_S^n(M))=0$, working from the last long exact sequence up to the
first long exact sequence, we obtain that
$\ext_{R_n}^{i+1}(S/\x,\tilde F_S^n(M))=0$. Similarly, one can
also obtain that $\ext_{R_n}^{i-1}(S/\x,\tilde F_S^n(M))=0$ if
$i>1$. Repeating these processes, we get
\[\ext_{R_n}^i(S/\x,\tilde F_S^n(M))=0\] for all $i>0$.

\medskip

Then, we show $\ext_R^i({}^{f^n}\!\! R,M)=0$ implies
$\ext_R^i({}^{f^{n+1}}\!\! R,M)=0, \forall \ i>0$. We need the
following lemma:
\begin{lem}Let $A \to B$ be a ring homomorphism such that $B$
is a finitely generated free module over $A$. Let $M$, $N$ be $A$-modules.
Then, as $A$-modules
\[
\ext_B^i(\hom_A(B,M),\hom_A(B,N)) \iso \hom_A(B,\ext_A^i(M,N)).
\]
\end{lem}
\begin{proof}[Proof of lemma]
Note that, as $A$-modules
\begin{align*}
\hom_B(\hom_A(B,M),\hom_A(B,N)) &\iso \hom_A(\hom_A(B,M)\otimes_B
B,N) \text{\ (adjointness)}\\
&\iso \hom_A(\hom_A(B,M),N) \\
&\iso \hom_A(M\otimes_A B, N) \text{\ \ \ \ \ \ \ \ (since $B$ is $A$-free)}\\
&\iso \hom_A(B, \hom_A(M,N)).
\end{align*}
This proves the case when $i=0$. For $i>0$, Let $J^\bullet$ be an
injective resolution of $N$. Then $\hom_A(B,J^\bullet)$ is an
injective resolution of $\hom_A(B,N)$ over $B$ since $B$ is free over $A$.
So, as $A$-modules
\begin{align*}
\ext_B^i(\hom_A(B,M),\hom_A(B,N)) &=
H^i(\hom_B(\hom_A(B,M),\hom_A(B,J^\bullet)))\\
&\iso H^i(\hom_A(B, \hom_A(M,J^\bullet)))\\
&\iso \hom_A(B,H^i(\hom_A(M,J^\bullet)))\text{\ ($B$ is $A$-free)}\\
&= \hom_A(B,\ext_A^i(M,N)).
\end{align*}
\end{proof}

\medskip

Now suppose for some fixed $i>0$ (therefore for all $i>0$ by the
first part of the proof)
\[
\ext_R^i({}^{f^n}\!\! R,M)=0,
\]
i.e.
\[\ext _{R_n}^i(S/\x, \tilde F_S^n(M))=0.\]

This implies that
\[
\hom_{R_n}({}^{\hat{f}}\! R_{n+1}, \ext_{R_n}^i(S/ \x,\tilde
F_S^n(M)))=0,\] where $\hat{f}$ denotes $f_S \otimes_S R_n:R_n \to
R_{n+1}$ (base change). Since ${}^{\hat{f}}\! R_{n+1}$ is free
over $R_n$, by the lemma, we have
\[
\ext_{R_{n+1}}^i(\hom_{R_n}({}^{\hat f}\!
R_{n+1},S/\x),\hom_{R_n}({}^{\hat f}\! R_{n+1},\tilde
F_S^n(M)))=0,
\]
i.e.
\[
\ext_{R_{n+1}}^i(\tilde F_S(S/\x), \tilde F_S^{n+1}(M))=0. \tag{3}
\]

\medskip

We claim that, as $S/\x^p$-modules, \[\tilde F_S(S/\x) \iso
S/\x^p.\]

To see this, first notice that by adjointness
\begin{align*}
         \tilde F_S(S/\x)&=\hom_S({}^{f}\! S, S/ \x)\\
                         &\iso \hom_{S/\x}(S/\x \otimes_S {}^{f}\!
                         S, S/\x)\\
                         &=\hom_{S/\x}({}^{\tilde f}(\! S/\x^p), S/\x),
\end{align*}
where $\tilde f$ is the map $\tilde f^n$ in the factorization (1)
with $n=1$. Then observe that since $S/\x$ and $S/\x^p$ are
complete intersections and $S/\x^p$ is a finitely generated
$S/\x$-module via $\tilde f$, $\hom_{S/\x}({}^{\tilde f}(\!
S/\x^p), S/\x)$ is a canonical module for $S/\x^p$ and hence is
isomorphic to $S/\x^p$ as an $S/\x^p$-module.

It follows that $\tilde F_S(S/\x)$ is isomorphic to $S/\x^p$ as an
$R_{n+1}$-module via the natural map $R_{n+1}\onto R_1$ (i.e. the
natural map $S/\x^{p^{n+1}} \onto S/\x^p$).

Replace $\tilde F_S(S/\x)$ in (3) by $S/\x^p$, we obtain
\[
\ext_{R_{n+1}}^i(S/{\x^p}, \tilde F_S^{n+1}(M))=0.\tag{4}
\]
Consider the following short exact sequences
\[
0 \to S/(x_1, x_2^p, \dots , x_r^p)
\overset{\lambda=x_1^p}{\longto} S/(x_1^{p+1}, x_2^p, \dots,
x_r^p) \to S/\x^p \to 0 \tag{5} \] and
\[
0 \to S/\x^p \overset{x_1}\to S/(x_1^{p+1}, x_2^p, \dots, x_r^p)
\overset{\mu}{\to} S/(x_1, x_2^p, \dots , x_r^p)\to 0. \tag{6}
\]
Apply $\hom_{R_{n+1}}(-,\tilde F_S^{n+1}(M))$ to (5). From the
associated long exact sequence and (4), we obtain an isomorphism
(induced by $\lambda$)
\[
\ext_{R_{n+1}}^i(S/(x_1,x_2^p,\dots,x_r^p),\tilde F_S^{n+1}(M))
\tag{7}\]
\[\iso \ext_{R_{n+1}}^i(S/(x_1^{p+1},x_2^p,\dots,x_r^p),\tilde
F_S^{n+1}(M)).
\]
Similarly, apply $\hom_{R_{n+1}}(-,\tilde F_S^{n+1}(M))$ to (6),
we get a map (induced by $\mu$)
\[
\ext_{R_{n+1}}^i(S/(x_1,x_2^p,\dots,x_r^p),\tilde F_S^{n+1}(M))
\tag {8}\]
\[\onto \ext_{R_{n+1}}^i(S/(x_1^{p+1},x_2^p,\dots,x_r^p),\tilde
F_S^{n+1}(M)),
\]
which is an isomorphism for $i \ge 2$ and a surjection for $i=1$.

The composition of maps $\lambda$ in (5) and $\mu$ in (6) is the
following multiplication
\[
S/(x_1,x_2^p,\dots,x_r^p) \overset {x_1^p} {\longrightarrow}
S/(x_1,x_2^p,\dots,x_r^p),
\]
which is a 0-map. By (7) and (8), it induces
\[
\ext_{R_{n+1}}^i(S/(x_1,x_2^p,\dots,x_r^p),\tilde F_S^{n+1}(M))
\overset {x_1^p}{\onto}
\ext_{R_{n+1}}^i(S/(x_1,x_2^p,\dots,x_r^p),\tilde F_S^{n+1}(M)).
\]
which is an isomorphism for $i\geq 2$ and an surjection for $i=1$.
Hence
\[
\ext_{R_{n+1}}^i(S/(x_1,x_2^p,\dots,x_r^p),\tilde F_S^{n+1}(M))=0.
\]
Repeating this process $(r-1)$ times, we obtain
\[
\ext_{R_{n+1}}^i(S/\x,\tilde F_S^{n+1}(M))=0,
\]
i.e.
\[
\ext_R^i({}^{f^{n+1}}\!\! R,M)=0
\]
for all $i>0$.

Therefore by induction, $\ext_R^i({}^{f^k}\!\! R,M)=0$ for all $k
\geq n$ and all $i>0$.

Finally, the assertion in the theorem follows immediately from the
following result due to Herzog ([H] Theorem 5.2): If $M$ is a
finitely generated module over $R$, then $M$ is of finite
injective dimension if and only if $\ext_R^i({}^{f^n}\!\! R,M)=0$
for all $i>0$ and infinitely many $n$.
\end{proof}

\begin{remark}The above theorem actually characterizes finitely generated
modules with finite projective dimension also since it is well
known that over Gorenstein ring, a finitely generated module is of
finite projective dimension if and only if it is of finite
injective dimension. To sum up, if $M$ is a finitely generated
module over a local complete intersection $R$, then the following
are
equivalent:\\
\begin{enumerate}
 \item [(i)] $\projdim_R M<\infty,$\\
 \item [(ii)] $\tor_i^R(M,{}^{f^n}\!\! R)=0$ for some fixed (all) $i,n>0$,\\
 \item [(iii)] $\ext_R^i({}^{f^n}\!\! R,M)=0$ for some fixed (all) $i,n>0$,\\
 \item [(iv)] $\injdim_R M<\infty$.\\
\end{enumerate}
\end{remark}
\begin{ac}I would like to thank professor Sankar Dutta for his
suggestions and comments. I also would like to thank the referee
for suggestions improving the presentation of this note.
\end{ac}

\bibliographystyle{amsalpha}

\end{document}